\documentclass{amsart}
\usepackage{amssymb, amsmath, amsfonts, amscd}

\input xypic

\theoremstyle{plain}
\newtheorem{theorem}{Theorem}[section]
\newtheorem{corollary}[theorem]{Corollary}
\newtheorem{lemma}[theorem]{Lemma}

\theoremstyle{definition}
\newtheorem{definition}[theorem]{Definition}

\theoremstyle{remark}

\numberwithin{equation}{theorem}

\newcommand{\I}{\mathcal{I}}

\newcommand{\E}{\mathcal{E}}
\renewcommand{\O}{\mathcal{O} }

\renewcommand{\P}{\mathbf{P} }

\renewcommand{\Pr}{\mathcal{P} }

\newcommand{\mf}[1]{\mathfrak{#1} }

\newcommand{\Spec}{\operatorname{Spec} }

\renewcommand{\H}{\operatorname{H} }

\renewcommand{\dim}{\operatorname{dim}}
\newcommand{\e}{\overline{e}}
\begin{document}

\title{A note on principal parts on projective space and linear representations}
\author{Helge Maakestad }
\address{Department of Mathematics, Bar-Ilan University,
 Ramat Gan, Israel }
\email{makesth@macs.biu.ac.il }
\thanks{Partially supported by the Emmy Noether Research Institute for Mathematics,
the Minerva Foundation of Germany, the Excellency Center ``Group Theoretic Methods in the
study of Algebraic Varieties'' of the Israel Science foundation and the EAGER foundation
(EU network, HPRN-CT-2000-00099)} 
\date{2002}
\begin{abstract} Let $H$ be a closed subgroup of a linear algebraic group 
$G$ defined over a field of characteristic zero. There is an equivalence of 
categories between the category of linear finite-dimensional representations of $H$,
and the category of finite rank $G$-homogeneous vector bundles on $G/H$. In this 
paper we will study this correspondence for the sheaves of principal parts on projective
space, and we describe the representation corresponding to the principal parts
of a line bundle on projective space.
\end{abstract}
\maketitle

\section{Introduction}

In this note we will study the vector bundles of principal parts 
$\Pr^k(\O (n))$ of a line bundle on projective space over a field $F$ of characteristic zero
from a representation theoretic point of view. We consider projective $N$-space
as a quotient $SL(V)/P$, where $V$ is an $N+1$-dimensional vectorspace over $F$,
and $P$ is the subgroup of $SL(V)$ stabilizing a line $L$ in $V$. There is 
an equivalence of categories between the category of finite rank 
$SL(V)$-homogeneous vector bundles on $SL(V)/P=\P(V^*)$ and the category of linear 
finite-dimensional representations of $P$. The principal parts $\Pr^k(\O(n))$
are $SL(V)$-homogeneous vector bundles on $\P(V^*)$, and the novelty of 
this note is that
we describe the $P$-representation corresponding to the principal parts.
The main result is Theorem \ref{maintheorem}. The Theorem says the following: 
Let $L^*$ be the dual of the $P$-module $L$. Then for all $1\leq k <n$, 
the $P$-representation
corresponding to $\Pr^k(\O(n))$ is $S^{n-k}(L^*)\otimes S^k(V^*)$. 
As a corollary, we obtain the splitting-type of $\Pr^k(\O(n))$ on $\P(V^*)$ 
for all
$1\leq k < n$, and recover results obtained in \cite{MAA},\cite{PERK2}, 
\cite{PIE} and \cite{dirocco}. 


\section{Principal parts on projective space}
 
In this section we give the representation corresponding to $\Pr^k(\O(n))$
on $\P(V^*)$ for all $1\leq k < n$, where $V$ is an $F$-vectorspace of dimension
$N+1$ and $F$ is a field of characteristic 0. A variety is an integral scheme of finite type 
over $F$. We will consider closed points when we talk about points of a scheme.
Let $V$ be a finite dimensional vectorspace over $F$. We let $GL(V)$ denote the group
of all invertible linear transformations of $V$. It is an algebraic group in the sense
of \cite{BOREL}, chapt. 1. A \emph{linear algebraic group} is a closed subgroup of
$GL(V)$. Let $SL(V)$ be the linear algebraic group of linear transformations of $V$ with
determinant $1$. Let $L$ in $V$ be a line, and
$P$ the closed  subgroup of $SL(V)$ stabilizing $L$. Then the quotient $SL(V)/P$ (which exists
by \cite{BOREL}, Theorem 6.8) is isomorphic
to $\P(V^*)$, the projective space of lines in $V$ (see \cite{AKH}, Section 4.2. This works 
over any field, not only the complex numbers). There exists a natural left
$SL(V)$-action on $\P(V^*)$, making it into a \emph{homogeneous space} for $SL(V)$. 
And by \cite{AKH}, Chapt. 4, there exists an equivalence of categories between 
the category of finite rank homogeneous
vector bundles on $SL(V)/P$ and the category of linear finite-dimensional representations
of $P$, and under this correspondence the dimension of the representation gives the rank
of the corresponding vector bundle. Hence any character of $P$ gives a homogeneous 
line bundle on $SL(V)/P$. The line $L$ corresponds to  
a character of $P$, and the bundle corresponding to the dual line $L^*$, is the line bundle
$\O(1)$ on $\P(V^*)$ (see \cite{AKH}, Section 4.2). It is also a standard fact that any linear 
finite-dimensional
representation $\rho$ of $P$ lifting to a representation $\tilde{\rho}$ of $SL(V)$ corresponds
to a trivial abstract vector bundle on $\P(V^*)$. There exists on any scheme an equivalence
of categories between the category of locally free finite rank sheaves and the category
of finite rank vectorbundles, hence we will use these two notions interchangably.

Pick a basis
$e_0,\dots ,e_N$ for $V$. Let $x_0,\dots ,x_N$ be the dual basis, and let $L$ be the 
line spanned 
by $e_0$. Having chosen a basis for $V$ it follows that
$SL(V)$ may be identified with the group of square rank $N+1$ matrices with 
determinant equal to 1. The group $P$ may be identified with the subgroup of $SL(V)$ 
consisting of matrices $g$ of the form
\[ g=
 \begin{pmatrix} a & * &    \cdots & * \\
                0 & a_{11} & \cdots & a_{1n} \\
                \vdots & \vdots   & \vdots & \vdots \\
                0 & a_{n1} & \cdots & a_{nn} \end{pmatrix} .\]
The one dimensional representation $\chi:P\rightarrow GL(S^n(L^*))$ corresponding 
to the line bundle $\O(n)$  is given by $\chi(g)=a^{-n}$.

 Let $X$ be a smooth variety of dimension $d$ and consider the diagonal $\Delta$ in $X\times X$. 
Let $\I$ be the sheaf of ideals of $\O_{X\times X}$ defining the diagonal $\Delta$,
and put $\O_{\Delta ^k}$ to be $\O_{X\times X}/\I^{k+1}$. 

\begin{definition} Let $p,q$ be the projetion maps from $X\times X$ to $X$, and let
$\E$ be an $\O_X$-module. Define $\Pr ^k(\E)=p_*(\O_{\Delta^k}\otimes q^*\E)$
to be the \emph{$k$'th order principal parts} of the module $\E$. We put
$\Pr^k(\O_X)=\Pr^k$. 
\end{definition}

Note that by \cite{MAA}, if the rank of $\E$ is $e$,
$\Pr^k(\E)$ is a vector bundle of rank $e\binom{d+k}{d}$ on $X$. 
Assume that $G$ is an 
algebraic group, and that $X$ is a homogeneous space for $G$. Assume furthermore
that $\E$ is a $G$-homogeneous vector bundle on $X$, then it follows that $\Pr^k(\E)$
is again a $G$-homogeneous vector bundle on $X$. Consider the line bundle $\O(n)$ on
$\P(V^*)=SL(V)/P$, then $\O(n)$ is an $SL(V)$-homogeneous line bundle on $\P(V^*)$
for all $n$, and we may consider the $SL(V)$-homogeneous vector bundle $\Pr^k(\O(n))$.
We want to compute the representation $\rho$ of $P$ corresponding
to the homogeneous vector bundle $\Pr^k(\O(n))$ for all $1\leq k <n$ on $\P(V^*)$. 
Let in the following $X=\P(V^*)$ and consider the projection maps $p,q$ 
from $X\times X$ to $X$.
Let $\I$ in $\O_{X\times X}$ be the ideal of the diagonal. We have an
exact sequence of $SL(V)$-homogeneous vector bundles on $X\times X$:

\begin{equation} \label{seq1}
\qquad 0 \rightarrow \I^{k+1} \rightarrow \O_{X\times X}
\rightarrow \O_{\Delta ^k} \rightarrow 0 .
\end{equation}
Apply the functor $p_*(-\otimes q^*\O(n))$ to the sequence \ref{seq1} 
to get a long exact sequence
\begin{equation} \label{seq2}
  \qquad0 \rightarrow p_*(\I^{k+1}\otimes q^*\O(n))\rightarrow 
p_*q^*\O(n)\rightarrow \Pr^k(\O(n)) 
\end{equation}
\[ \rightarrow R^1p_*(\I^{k+1}\otimes q^*\O(n))\rightarrow R^1p_*q^*\O(n)\rightarrow 
R^1p_*(\O_{\Delta^k}\otimes q^*\O(n)) \rightarrow \cdots \]
of vector bundles. The sequence \ref{seq2} is a sequence of vector bundles 
because all sheaves in the sequence are coherent, and it is a standard 
fact that a homogeneous coherent sheaf on a homogeneous space is locally free.
Since the sequence  \ref{seq2} is a sequence of vector bundles, we 
get an exact sequence of $P$-representations when we pass to the fiber at $\overline{e}$. 


Consider the diagram
\[ \diagram \Spec(\kappa(\overline{e}))\times X\dto^\pi  \rto^j & 
X\times X \dto^p \\
 \Spec(\kappa(\overline{e})) \rto^i & X \enddiagram
,\]
then by \cite{HH}, Chapt. III, Sect.12 we get maps
\[\phi^i: R^ip_*(\I^{k+1}\otimes q^*\O(n))(\overline{e}) \rightarrow R^i\pi_*(j^*(\I^{k+1}\otimes q^*\O(n))) \]
of $\O_X$-modules. Put for any $\O_{X\times X}$-module $\E$ 
\[ h^i(y,E)=\dim_{\kappa (y)}\H^i(X_y,\E_y) ,\]
where $X_y$ is the fiber $p^{-1}(y)$ and $\E_y$ is the restriction of $\E$ to
$X_y$.
We see that 
\[ h^i(y,\I^{k+1}\otimes q^*\O(n))=\dim_{\kappa (y)}
\H^i(X,\mf{m}_y^{k+1}\otimes \O(n)) \]
is a constant function of $y$ for $i=0,1,\dots$ for the following reason:
Consider the following commutative diagram
\[
\diagram \Spec(\kappa (y))\times X \dto^j \rto^{\tilde{g} } & \Spec(\kappa(gv))\times X
\dto^k \\
X \times X \rto^g & X \times X \enddiagram \]
where the action of $SL(V)$ on $X\times X$ is given by
$g(x,y)=(gx,gy)$. In general if $G\times Y\rightarrow^\sigma Y$  is an algebraic group acting 
on a scheme $Y$, and $\E$ is a $G$-linearized sheaf on $Y$, then there
exists an isomorphism $I:\sigma^*\E\rightarrow p^*\E$ where $p:G\times Y\rightarrow Y$
is the projection map. It follows that for all $g\in G$ we get an isomorphism
$g^*\E\cong \E$ of sheaves.
Then since $\O(n)$ and $\I^{k+1}\otimes q^*\O(n)$
are $SL(V)$-homogeneous sheaves, we have an isomorphism
\[ \tilde{g}^*(\mf{m}_{gy}^{k+1}\otimes \O(n))=j^*g^*(\I^{k+1}\otimes 
q^*\O(n))=\mf{m}_y^{k+1}\otimes \O(n) \]
hence since $\tilde{g}$ is an isomorphism, 
we see that  we have an isomorphism
\[ \mf{m}_y^{k+1}\otimes \O(n)\cong \mf{m}_{gy}^{k+1}\otimes \O(n) \]
of sheaves for all $g$ in $SL(V)$. It follows that
\[  \dim_{ \kappa (y) }\H^i(X,\mf{m}_y^{k+1}\otimes \O(n) )= 
\dim_{ \kappa (gy)}\H^i(X,\mf{m}_{gy}^{k+1},\O(n) ) \]
for all $g$ in $SL(V)$,
hence  by \cite{HH}, chapter III Corr. 12.9, it follows that the maps $\phi^i$ 
are isomorphisms for $i=0,1,\dots$. Here $\mf{m}_y$ is the sheaf of ideals
corresponding to the point $y$ in $X$.
We get an exact sequence 
\begin{equation} \label{seq3}
\qquad 0 \rightarrow \H^0(X,\O(n)\otimes \mathfrak{m}_{\e}^{k+1})\rightarrow 
\H^0(X,\O(n))\rightarrow
\Pr^k(\O(n))(\overline{e}) 
\end{equation}
\[\rightarrow \H^1(X, \O(n)\otimes \mathfrak{m}_{\e}^{k+1}) \rightarrow \H^1(X, \O(n)) 
\rightarrow \cdots \]
of $P$-representations.

\begin{lemma} \label{lemma1} For all $1\leq k <n$ we have that $\H^1(X,\O(n)\otimes 
\mathfrak{m}_{\e}^{k+1})=0$.
\end{lemma}
\begin{proof} Consider the exact sequence \ref{seq3}. We prove that
\[\dim _F \H^0(X,\O(n))-\dim_F \H^0(X,\O(n)\otimes \mathfrak{m}_{\e}^{k+1}) = 
\dim _F \Pr^k(\O(n))(\overline{e}),\]
and then the result follows by counting dimensions.
We have that $\H^0(X,\O(n))$ equals $S^n(V^*)$, where $V^*$ is the $F$-vectorspace
on the basis $x_0,\dots ,x_N$. We also see that $\H^0(X,\O(n)\otimes 
\mathfrak{m}_{\e}^{k+1})$
equals $m^{k+1}S^{n-(k+1)}(V^*)$ considered as a subvectorspace of $S^n(V^*)$.
 Here $m$ is the $F$-vectorspace on the basis $x_1,\dots ,x_N$
and $m^{k+1}S^{n-(k+1)}(V^*)$ is the image of the natural map
\[ S^{k+1}(m)\otimes S^{n-(k+1)}(V^*)\rightarrow S^n(V^*).\]
 Write $V^*$ as
the direct sum $Fx_0\oplus m$. Then it follows that
\[ m^{k+1}S^{n-(k+1)}(V^*)=x_0^{n-(k+1)} m^{k+1}\oplus \cdots \oplus 
x_0 m^{n-1}\oplus m^n ,\]
hence we see that the dimension of $m^{k+1}S^{n-(k+1)}(V^*)$ equals
$\sum_{i=k+1}^n \binom{i+N-1}{N-1}$. We also see that the dimension
of $S^n(V^*)$ equals $\sum_{i=0}^n \binom{i+N-1}{N-1}$, and it follows
that 
\[ \dim _F S^n(V^*)-\dim _F m^{k+1}S^{n-(k+1)}(V^*) = \sum_{i=0}^k 
\binom{i+N-1}{N-1}.\]
It follows that  
\[ \sum_{i=0}^k \binom{i+N-1}{N-1}= \binom{k+N}{N}=\dim _F \Pr^k(\O(n))(\overline{e}),\]
and we have proved that
\[\dim _F\H^0(X,\O(n))-\dim _F\H^0(X,\O(n)\otimes \mf{m}_{\e}^{k+1})=
\dim _F \Pr^k(\O(n))(\overline{e}),\] 
and the result follows from the fact that the sequence \ref{seq3} is exact and that
$\H^1(X,\O(n))=0$ for $n\geq 1$.
\end{proof}

Note that by Lemma \ref{lemma1} and the sequence \ref{seq3}, there exists for 
all $1\leq k <n$ an exact sequence
of $P$-representations
\[ 0 \rightarrow \H^0(X,\O(n)\otimes \mf{m}_{\e}^{k+1})\rightarrow 
\H^0(X,\O(n))\rightarrow
\Pr^k(\O(n))(\overline{e}) \rightarrow 0 .\]
Since the representation 
$\H^0(X, \O(n)\otimes  \mf{m}_{\e}^{k+1})$ equals $m^{k+1}S^{n-(k+1)}(V^*)$ 
as subrepresentation
of $\H^0(X, \O(n))=S^n(V^*)$, it follows that we have an exact sequence of $P$-representations
\begin{equation} \label{seq4}
0\rightarrow m^{k+1}S^{n-(k+1)}(V^*)\rightarrow S^n(V^*)\rightarrow
\Pr^k(\O(n))(\overline{e})\rightarrow 0.
\end{equation} 

From the exact sequence

\[ 0 \rightarrow m \rightarrow V^* \rightarrow V^*/ m \rightarrow 0 ,\]
where $m$ is the $F$-vectorspace on $x_1,\dots ,x_N$,
we see that the representation $V^*/m$ is the representation  corresponding to the
module  $L^*$ of $P$, giving the line bundle $\O(1)$ on $X=SL(V)/P$. 

\begin{lemma} \label{lemma2} For all $1\leq k <n$ there exists a surjective 
map of $P$-representations
\[ \phi:S^n(V^*)\rightarrow S^{n-k}(L^*)\otimes S^k(V^*) .\]
\end{lemma}
\begin{proof} Recall that we have chosen a basis $e_0,\dots ,e_N$ for $V$, with the property that
$\overline{x}_0$ is a basis for $L^*$. The $P$-representation $m$ with basis 
$x_1,\dots ,x_N$ gives an exact sequence
\[ 0\rightarrow m \rightarrow V^* \rightarrow L^* \rightarrow 0 \]
of $P$-representations. Define a map 
\[ \phi:S^n(V^*)\rightarrow S^{n-k}(L^*)\otimes S^k(V^*) \]
as follows: $\phi(f)=\overline{x}_0^{n-k}\otimes \partial_0^{n-k}(f)$, where $\partial_0^{n-k}$ 
is $n-k$ times partial derivative with respect to the $x_0$-variable. 
Let $g$ be an element of $P$. Then by
induction on the degree of the differential-operator $\partial_0^{n-k}$ and applying
the chain-rule for derivation, it follows that 
\[\phi(gf)=\overline{x}_0^{n-k}\otimes \partial_0^{n-k}(gf)=
\overline{x}_0^{n-k}\otimes a^{-(n-k)}g(\partial_0^{n-k}f)= \]
\[ a^{-(n-k)}\overline{x}_0^{n-k}\otimes g(\partial_0^{n-k}f)=g(\overline{x}_0^{-(n-k)}\otimes
\partial_0^{n-k}f)=g\phi(f) ,\]
and we see that $\phi$ is $P$-linear. It is clearly surjective, and the lemma follows.
\end{proof}

\begin{theorem} \label{maintheorem} For all $1\leq k <n$, 
the representation corresponding to $\Pr^k(\O(n))$ is 
$S^{n-k}(L^*)\otimes S^k(V^*)$. 
\end{theorem}
\begin{proof} By Lemma \ref{lemma2} there exists a surjective map 
of $P$-representations 
\[ \phi :S^n(V^*)\rightarrow S^{n-k}(L^*)\otimes S^k(V^*).\]
We claim that 
$m^{k+1}S^{n-(k+1)}(V^*)$ equals $ker \phi$: We first prove the inclusion
\[ m^{k+1}S^{n-(k+1)}(V^*)\subseteq ker \phi. \]
Pick a monomial $x_0^{p_0}x_1^{p_1}
\cdots x_N^{p_N}$ in $m^{k+1}S^{n-(k+1)}(V^*)$, hence 
$p_0+\cdots +p_N=n$ and $p_0 <n-k$. These monomials form a basis for 
$ m^{k+1}S^{n-(k+1)}(V^*)$. We see that $\partial_0^{n-k}(x_0^{p_0}\cdots 
x_N^{p_N})$ is zero, hence since $\phi$ is a linear map, it follows that 
we have an inclusion 
\[ m^{k+1}S^{n-(k+1)}(V^*)\subseteq ker \phi \] 
of vectorspaces.
The reverse inclusion follows from counting dimensions and the fact that $\phi$ is
surjective: We have that
\[ \dim _F ker\phi=\dim _F S^n(V^*)-\dim _F S^{n-k}(V^*)\otimes S^k(V^*) \]
\[=\sum_{i=0}^n\binom{i+N-1}{N-1}-\sum_{i=0}^k\binom{i+N-1}{N-1}=\sum_{i=k+1}^n\binom{i+N-1}{N-1},\]
and we see that $\dim _F ker\phi= \dim _F m^{k+1}S^{n-(k+1)}(V^*)$ and it follows that
\[m^{k+1}S^{n-(k+1)}(V^*)=ker\phi,\] 
hence we have an exact sequence of $P$-representations
\[ 0 \rightarrow m^{k+1}S^{n-(k+1)}(V^*)\rightarrow S^n(V^*)
\rightarrow S^{n-k}(L^*)\otimes S^k(V^*) \rightarrow 0 .\]
Using sequence \ref{seq4} we get isomorphisms
\[ \Pr^k(\O(n))(\overline{e})\cong \H^0(X, \O(n))/\H^0(X,\O(n)\otimes 
\mf{m}_{\e}^{k+1})
\cong  \]
\[ S^n(V^*)/m^{k+1}S^{n-(k+1)}(V^*)\cong S^{n-k}(V^*)\otimes S^n(V^*) ,\]
and it follows that $\Pr^k(\O(n))(\overline{e})$ is isomorphic to
$S^{n-k}(L^*)\otimes S^k(V^*)$ as representation.
\end{proof}

Note that the result in Theorem \ref{maintheorem} is true if $char(F)>n$.

\begin{corollary} For all $1\leq k<n$, $\Pr^k(\O(n))$ splits as abstract
vector bundle as $\oplus^{\binom{N+k}{N} }\O(n-k)$.
\end{corollary}
\begin{proof} Since $S^k(V^*)$ corresponds to the trivial rank
$\binom{N+k}{N}$ abstract vector bundle on $\P(V^*)$, and $S^{n-k}(L^*)$ 
corresponds 
to the line bundle $\O(n-k)$, the assertion is proved.
\end{proof}

We see that we recover results on the splitting-type of the principal parts
obtained in \cite{MAA}, \cite{PERK2},\cite{PIE} and \cite{dirocco}.

\textbf{Acknowledgements}. I would like to thank Michel Brion for an 
invitation to spend spring 2001 at the Institut Joseph Fourier.
Thanks for suggesting the problem, sharing ideas
 and for answering numerous questions on homogeneous spaces and representation
theory.  Thanks also to Laurent Manivel for stimulating 
discussions on the subject. 
I would also like to thank Mina Teicher for an invitation to the 
Bar-Ilan University, where parts of this work was done and this note was written.
Finally, thanks to Dan Laksov for comments.


\begin{thebibliography}{4}



\bibitem{AKH} D. N. Akhiezer, Lie group actions in complex analysis,
Aspects of Mathematics, \emph{Vieweg} (1995)

\bibitem{BOREL} A. Borel, Linear algebraic groups, Graduate Texts in
Mathematics no. 126,
\emph{Springer Verlag} (1991)




\bibitem{JAN} J. C. Jantzen, Representations of algebraic groups,
Pure and Applied Mathematics no. 131, \emph{Academic Press} (1987)

\bibitem{HH} R. Hartshorne, Algebraic geometry, Graduate Texts in 
Mathematics no. 52, \emph{Springer Verlag} (1977)

\bibitem{MAA} H. Maakestad, Modules of principal 
parts on the projective line, \emph{preprint math.AG/0111149}



\bibitem{PERK2} D. Perkinson, Principal parts of linebundles on
    toric varieties, \emph{Compositio Math.} 104, 27-39, (1996)

\bibitem{PIE} R. Piene, G. Sacchiero, Duality for rational
    normal scrolls, \emph{Comm. in Alg.} 12 (9), 1041-1066, (1984) 


\bibitem{dirocco} S. di Rocco, A. J. Sommese, Line bundles for
which a projectivized jet bundle is a product,
\emph{Proc. Amer. Math. Soc.} 129, (2001), no.6, 1659--1663 

\end{thebibliography}
\end{document}